\documentclass[onecolumn]{autart}
\usepackage{}    
\usepackage[colorlinks,linkcolor=black,anchorcolor=red,citecolor=red]{hyperref}%

\usepackage{graphicx}           
\usepackage{dcolumn}            
\usepackage{bm}                 
\usepackage{bbm}
\usepackage{graphics}           
\usepackage{epsfig}             
\usepackage{times}              
\usepackage[fleqn]{amsmath}
\usepackage{amssymb}
\usepackage{extarrows}
\usepackage{amsfonts}
\usepackage{mathrsfs}
\usepackage{fancyhdr}
\usepackage{color}
\usepackage{subfigure}
\usepackage{enumerate}
\usepackage [latin1]{inputenc} %
\allowdisplaybreaks[4] 

\newtheorem{theorem}{Theorem}[section] %
\newtheorem{lemma}{Lemma}[section]

\newtheorem{proposition}[theorem]{Proposition}

\newtheorem{remark}{Remark}[section]

\DeclareMathOperator{\esssup}{ess\,sup}
\begin{document}

\begin{frontmatter}

\title{Exponential stabilization and  continuous dependence  of solutions on  initial data in different norms for   space-time-varying linear  parabolic PDEs}

\author{Qiaoling Chen$^{1}$}\ead{cql20202016@my.swjtu.edu.cn},
\author{Jun Zheng$^{1,2}$}\ead{zhengjun2014@aliyun.com},
\author{Guchuan Zhu$^{2}$}\ead{guchuan.zhu@polymtl.ca}

\address{$^{1}${School of Mathematics, Southwest Jiaotong University,
        Chengdu 611756, Sichuan, China}\\
        $^{2}$Department of Electrical Engineering, Polytechnique Montr\'{e}al, P.O. Box 6079, Station Centre-Ville, Montreal, QC, Canada H3T 1J4}

\begin{abstract}
 For an arbitrary parameter $p\in [1,+\infty]$, we   consider the problem of exponential stabilization in the spatial $L^{p}$-norm, {and $W^{1,p}$-norm, respectively,} for a class of anti-stable  linear  parabolic PDEs with    space-time-varying coefficients  in the absence of a Gevrey-like condition, which is often imposed  on  time-varying  coefficients of PDEs and used to guarantee the existence of  smooth (w.r.t. the time variable) kernel functions in the literature.  Then, based on the obtained exponential stabilities, we show  that the solution of the considered system depends continuously
   on the  $L^{p}$-norm,  and $W^{1,p}$-norm, respectively, of the initial data.     In order to obtain   time-independent  (and thus  sufficiently smooth) kernel functions  without   a Gevrey-like condition and  deal with singularities arising in  the case of $p\in[1,2)$, we apply   a combinatorial method, i.e., the combination  of backstepping and approximation  of Lyapunov functionals (ALFs), to stabilize the considered system  and establish the continuous dependence of  solutions on   initial data in different norms.

\end{abstract}
\end{frontmatter}
 \section{Introduction}

In the last decades,   the backstepping method has   proved to be  a powerful tool for the boundary control of distributed parameter systems governed by  $1$-D partial differential equations (PDEs); see, e.g., \cite{krstic2008boundary}.
The main idea of the method is to transform the  anti-stable PDEs   into  stable target systems by using a boundary feedback control  and  a change of variable of the PDEs, i.e., two integral transformations, so as to eliminate the anti-stable {term} of the considered PDEs. The process of conducting stability analysis via this method mainly consists of two steps: (i) proving the existence of the needed transformation  and its inverse in a certain function space,   which are often  accordant with  the existence and regularity of two kernel functions, and can be proceeded in the standard way; and (ii) establishing the  stability  estimates for the  target systems, which may be  based on different methods of PDEs or functional analysis.

In the context of parabolic systems,   significant progress  has been obtained  in exponential    stabilization in   $L^2$-norm  or $H^{1}$-norm  for linear parabolic PDEs by using the method of backstepping; see, e.g., \cite{krstic2008boundary,boskovic2001boundary,elharfi2008explicit,Feng:2022,Karafyllis:2021,Kerschbaum:2020,Kerschbaum:2022,liu2003boundary,meurer2009tracking,si2018boundary,smyshlyaev2005control,Smyshlyaev:2004,balogh2004stability,smyshlyaev2005backstepping}. In general, when  transformations are applied,    the target systems are   expected to be in a simple form, whose stability can be established easily. For example,  the  target systems usually have the following form:
\begin{align}\label{equ.example1}
u_t(x,t)=u_{xx}(x,t)-\lambda  u(x,t),
\end{align}
where    $\lambda$ is an arbitrary constant; see, e.g., \cite{krstic2008boundary,boskovic2001boundary,elharfi2008explicit,Feng:2022,Karafyllis:2021,Kerschbaum:2020,Kerschbaum:2022,liu2003boundary,meurer2009tracking,si2018boundary,smyshlyaev2005control}. In particular, for $\lambda>0$, the stability in different norms of \eqref{equ.example1} has been well studied, and  is    used together with the inverse transformations to obtain  the stability in $L^2$-norm  or $H^{1}$-norm  for the original parabolic PDEs. Moreover, since the systems are linear, it is clear that the stability in $L^2$-norm  or $H^{1}$-norm implies   the  continuous dependence of    solutions    on the initial data in $L^2$-norm  or $H^{1}$-norm.

It should be mentioned that, for parabolic PDEs with  time-varying coefficients,   certain Gevrey-like condition, i.e., the time-varying coefficients of the reaction terms   are supposed to be in Gevrey class, is often imposed   to guarantee the existence of  sufficiently   smooth (w.r.t. the time variable) kernel functions when the method of {backstepping} is applied; see, e.g., \cite{Kerschbaum:2020,meurer2009tracking,si2018boundary,smyshlyaev2005control}. It is reasonable to believe that the {method}  of {backstepping} along with  the Gevrey-like condition  can be adopted to  the exponential     stabilization  in   $L^p$-norm   {or $W^{1,p}$-norm}, and to the establishment of  continuous dependence of    solutions    on the initial data in $L^p$-norm   {or $W^{1,p}$-norm},  for linear parabolic PDEs  with space-time-varying coefficients  whenever $p\in [1,+\infty]$. Nevertheless,  it  is challenging to obtain sufficiently smooth (w.r.t. the time variable) kernel functions via   backstepping if   no Gevrey-like condition is imposed on time-varying coefficients. Indeed, for a target system  chosen as \eqref{equ.example1}, the   equations of kernel functions  are relevant to the reaction coefficients of the original systems, and the proof of existence and higher regularity of   kernel functions depends  strongly on the application of the property of the reaction coefficients, namely, a Gevrey-like condition seems  to be needed; see, e.g., \cite{Kerschbaum:2020,meurer2009tracking,si2018boundary,smyshlyaev2005control}.

An alternative is to expect that the equations of kernel functions  are independent of the time variable,  which may extensively reduce  the complications of deriving kernel functions. In this case,  a target system  should have the following form:
\begin{align} \label{equ.example}
u_t(x,t)=u_{xx}(x,t)-\lambda(x,t) u(x,t),
\end{align}
    where $\lambda(x,t)>0$ depends on the reaction coefficients of the original system, and hence is space-time varying.  However, for   {space-time varying $\lambda(x,t)$}, it is   difficult to establish the stability and the continuous dependence of solutions on the initial data
   involved with the  spatial   $L^p$-norm {or $W^{1,p}$-norm} for system \eqref{equ.example} whenever $p\in [1,+\infty]$.

    Let  us {take an example by considering} system \eqref{equ.example}  {that is} defined over $(0,1)\times(0,+\infty)$,  and with  homogeneous Dirichlet boundary conditions and $L^p$-initial data $u_0$.
    It is well-known that the Lyapunov method is a well-suited tool for  stability  analysis of parabolic PDEs. In particular, for the linear system \eqref{equ.example}, the functional $$V(u):=\int_0^1|u(x,t)|^p\text{d}x $$ is an suitable Lyapunov candidate that can be used  to establish  the  exponential stability in  $L^p$-norm  having the form
    \begin{align} \label{Lp}
    \|u[t]\|_{L^p(0,1)}\leq C_1 e^{\sigma t}  \|u_0\|_{L^p(0,1)},\forall t\geq 0,
     \end{align}
   when $p\geq 2$, where $C_1>0$ and $\sigma>0$ are some constants; see, e.g., \cite{krstic2008boundary,Zheng:201802}.

    In the case of $p\in[1,2)$, it may be expected that  the  functional $V(u)$ could be  used to  establish the stability in  $L^p$-norm  of   system \eqref{equ.example}. However, as indicated in \cite{zheng2021approximations,Zheng:202112},  the choice of $V(u)$   may lead  to singularities when $p\in[1,2)$.    {Indeed,}
      \begin{enumerate}
\item[$\bullet$]for $p=1$, due to the fact that the function $g(s):=|s|$ is not differentiable  at $s=0$,   the functional
$
 {\frac{\text{d}}{\text{d}t}\int_0^1|u |\text{d}x}$ is singular;
\item[$\bullet$]for $1<p<2$, by integrating by parts, it holds (formally) that
\begin{align*}
 \frac{\text{d}}{\text{d}t}\int_0^1|u|^p\text{d}x=&p\int_0^1|u|^{p-1}\text{sgn}(u)u_t\text{d}x\\
=&p\int_0^1|u|^{p-1}\text{sgn}(u)(u_{xx} -\lambda u)\text{d}x\\
=&-p\int_0^1\frac{\text{d}}{\text{d}u}(|u|^{p-1}\text{sgn}(u))u_x^2\text{d}x - p\int_0^1\lambda|u|^{p}  \text{d}x.
\end{align*}
But $\frac{\text{d}}{\text{d}u}(|u|^{p-1}\text{sgn}(u))$ is singular at $u=0$, where $\text{sgn}(\cdot)$ denotes the standard sign function.
\end{enumerate}

      It may be also expected that the stability in $L^p$-norm could be established by using the Sobolev embedding $L^2(0,1)\hookrightarrow L^p(0,1)$.  For instance, for $p\in [1,2)$,     it follows from  $L^2(0,1)\hookrightarrow L^p(0,1)$ and \eqref{Lp} that
    \begin{align*}
    \|u[t]\|_{L^p(0,1)}\leq &C_2 \|u[t]\|_{L^2(0,1)} 
       \leq  C_1C_2 e^{\sigma t}  \|u_0\|_{L^2(0,1)},\forall t\geq 0,
     \end{align*}
     where $C_2>0$ is the Sobolev embedding constant\footnote{It can be also determined by using the H\"{o}lder's inequality.}. However,   the stability estimate  obtained in this way involves  the $L^2$-norm of the initial data,
        which does not imply the continuous dependence of solutions on the   $L^p$-norm of the initial data.

 It {is   worth} mentioning that, compared with the {upper} case of $p\in[2,+\infty]$,  the intermediate  case of $p\in (1,2)$ and the lower case of $p=1$    {make}  the problem  more meaningful.  For example,  many physical models lead to PDEs with $L^1$-data and should be estimated in $L^1$-norm, or   need  to  be shown that the solutions  are continuously dependent on the $L^1$-data; {see, e.g., \cite{Goudon:1998,Lewandowski:1997,Li:2001,Roubicek:1998}}. It is also worth noting that the  study on the continuous  dependence of solutions on   initial data for PDEs is  important not only for the theory of regularity   of evolution equations, but also for the analysis of stability and robustness of   PDE systems, e.g., the input-to-state stability (ISS) and  relative stability; see \cite{Zheng:2022} for instance, where by using the De Giorgi iteration  the continuous  dependence of solutions on  the $L^{\infty}$-norm of initial data was established for a class of nonlinear parabolic PDEs with Robin boundary conditions, and used to prove the ISS in $L^{\infty}$-norm for the considered systems and the stability in $L^{\infty}$-norm for a class of cascaded nonlinear
parabolic systems, respectively.

  Nevertheless, due to the  aforementioned difficulties, few literature concerns the stabilization and continuous dependence  of solutions on  initial data in the spatial $L^p$-norm ({or $W^{1,p}$-norm}) for    anti-stable parabolic PDEs with time-varying coefficients  when $p $ is restricted to be  in $[1,2)$  and smooth kernel functions are adopted  without  a Gevrey-like condition in backstepping. This is the main motivation of the work.


In this paper,  without any Gevrey-like condition and for arbitrary $p\in [1,+\infty]$, we  consider the problem of exponential stabilization and continuous dependence  of solutions on  initial data  {for a class of space-time-varying linear parabolic PDEs,  and establish for solutions the estimates  containing  the spatial $L^p$-norm,  and $W^{1,p}$-norm, of the initial data, respectively}. The method adopted in this paper is based on the combination of   backstepping  and  approximations of Lyapunov functionals (ALFs),  among which the latter were used in, e.g.,   \cite{zheng2021approximations,Zheng:202112}, to analyze stability of nonlinear  PDEs with integrable inputs.  More specifically,
\begin{enumerate}
 \item[$\bullet$]  in order to eliminate  the anti-stable term of the considered parabolic PDEs,   we design a boundary feedback  control    by using the {method of backstepping};
       \item[$\bullet$]  in the absence of a Gevrey-like condition, we  consider time independent kernel functions, and transform the original system   into certain target system  that has  space-time-varying coefficients;
             \item[$\bullet$] in order to {deal with} singularities that may arise in  the case of $p\in[1,2)$, we apply the {technique}  of  {ALFs} to analyze the stability and continuous dependence  of solutions on  initial data  in different norms for the target system  whenever $p\in[1,+\infty]$.
\end{enumerate}

The main contribution of this paper is to apply  such a combinatorial method to  establish various estimates of solutions for a class of anti-stable parabolic PDEs, and is two-fold:
 \begin{enumerate}
 \item[$\bullet$]   exponential stability and continuous dependence of solutions on initial data are established in different norms;

             \item[$\bullet$] kernel functions are allowed to be independent of the time variable, and therefore the Gevrey-like condition can be avoided.
\end{enumerate}

In the rest of the paper, we introduce first some basic
notations.  In Section \ref{Sec. 2}, we present the problem formulation,    main result  and technical line of the proof.   In Section \ref{Sec. 3}, by using  the {technique}  of {ALFs}, we analyze the   exponential  stability and continuous dependence  of solutions on  initial data in different norms for the target systems. In Section \ref{Sec. 4}, we present  a result of well-posedness and regularity of kernel functions. In Section \ref{Sec. 5}, we prove the main result stated in Section \ref{Sec. 2}. Some conclusions are given in Section \ref{Sec. 6}.

\paragraph*{Notation.}  {Let $\mathbb{N}_0$}, $\mathbb{N}$, $\mathbb{R}$, {$\mathbb{R}_{\geq 0}$} and $ {\mathbb{R}_{>0}}$ be the set of nonnegative integers,  positive integers, real numbers, {nonnegative real numbers,} and positive real numbers, respectively.

 For $q\in[1,+\infty]$,  the notations of $ L^q(0,1)$
 and  $ W^{1,q}(0,1)$  denote the standard Lebesgue spaces and Sobolev spaces with elements  defined over $(0,1)$  and equipped with the norm $ {\|v\|_{q}:=}\|v\|_{L^q(0,1)} $ and
${\|v\|_{1,q}}:=\|v\|_{W^{1,q}(0,1)} $, respectively; see, e.g., \cite{Evans:2010}.
 %

Let $L^{\infty}((0,1)\times {\mathbb{R}_{>0}}):=\{v:(0,1)\times {\mathbb{R}_{>0}} \rightarrow \mathbb{R}| ~v$ is measurable and satisfies $\esssup_{(y,t)\in(0,1)\times {\mathbb{R}_{>0}}}|v{(y,t)}|<+\infty\}$.

 For   $v: (0,1)\times {\mathbb{R}_{>0}}\rightarrow \mathbb{R}$, the notation $v[t]$ {(or $v[y]$)  denotes the profile at certain $t\in{\mathbb{R}_{>0}}$ (or $y\in(0,1)$), i.e., $ v[t] (y)=v(y,t)$ (or $ v[y] (t)=v(y,t)$) for all $y\in(0,1)$ (or $t\in\mathbb{R}_{>0}$)}. For $T\in{\mathbb{R}_{>0}}$ and  a Banach space $X$ equipped with the norm $\|\cdot\|_{X}$, let $C\left([0,T];X\right):=\{v:(0,1)\times[0,T]\rightarrow \mathbb{R} |~v[t]\in X $ and satisfies $\max_{t\in[0,T]}\|v[t]\|_{X }<+\infty \}$, which is a   Banach space equipped with the norm $\|v\|_{C\left([0,T];X\right)}:=\max_{t\in[0,T]}\|v[t]\|_{X}$.

For $Q ={\mathbb{R}_{\geq0}}$ (or $[0,1]$), let $C^1\left([0,1]\times Q\right):=\{v:[0,1]\times Q \rightarrow\mathbb{R}|~v$ has continuous derivatives up to order $1$$\}$.

 For $l\in{\mathbb{R}_{>0}}$ and $T \in{\mathbb{R}_{>0}}$, the notions of $\mathcal{H}^{2+l}\left([0,1]\right)$ and $\mathcal{H}^{2+l,1+\frac{l}{2}}\left([0,1]\times [0,T]\right)$ denotes H\"{o}lder spaces that are defined by \cite[Chapter \uppercase\expandafter{\romannumeral1}]{Ladyzhenskaya:1968}.

\section{Problem Setting,~~Main Result and Technical Line}\label{Sec. 2}
\subsection{Problem Setting and Main Result}
 Given certain initial data $w_0$ , we study the problem of stabilization and continuous dependence  of solutions on  initial data  in  different norms for     the  linear parabolic PDEs with   space and time dependent {coefficients}:
 \begin{subequations}\label{original system}
\begin{align}
	 w_{t}(x,t)=&w_{xx}(x,t)+ c(x,t)w(x,t)  +\int_0^xw(y,t)f(x,y)\text{d}y, (x,t)\in (0,1)\times {\mathbb{R}_{>0}},\label{1a}\\
	 w_x(0,t)=&0, t\in {\mathbb{R}_{>0}},\label{1b}\\
	 w_x(1,t)=&U(t),t\in {\mathbb{R}_{>0}}, \label{1c}\\
	 w(x,0)=&w_0(x),x\in(0,1),
\end{align}
\end{subequations}
where  $c: [0,1]\times {\mathbb{R}_{\geq0}} \rightarrow \mathbb{R} $ and $f:[0,1]\times[0,1]\rightarrow \mathbb{R} $     are given functions,   and  $U$ is the control input to be determined to stabilize the system.

For the system \eqref{original system}, the reactive coefficient $c$ often represents characteristic
quantities for heat exchange and reaction, {e.g.,  \cite{meurer2005feedforward,jensen1982bifurcation}.} Moreover, if $c$ is positive and {sufficiently large}, the open-loop system (i.e., $U(t)=0$) is unstable; see, e.g., \cite{krstic2008boundary}.

 Throughout this paper, we assume that
 \begin{align*}
 &{p\in  [1,+\infty)}\ \text{is\ an\ arbitrary\ constant,\ or}\ p=+\infty,\\
 &f\in  C^1 \left([0,1]\times[0,1]\right), c\in  C^1\left([0,1]\times {\mathbb{R}_{\geq0}}\right)\cap L^{\infty}((0,1)\times {\mathbb{R}_{>0}}).
  \end{align*}
Moreover, assume that $c$ has a form
\begin{align}
 c(x,t):=&c_1(x)+c_2(t)\label{c}
 \end{align}
with  some functions $c_1:[0,1] \rightarrow \mathbb{R}~\text{and}~c_2:   {\mathbb{R}_{\geq 0}}\rightarrow \mathbb{R}$.

 It is worth noting that stabilization of anti-stable parabolic PDEs   under such kind of structural condition \eqref{c} was firstly considered in \cite{smyshlyaev2005control}, where  boundary feedback controls were designed to exponentially stabilize   equation \eqref{1a}  (with $f(x,t)\equiv 0$)  in  $L^2$ and $W^{1,2}$ norms for   $ c(x,t)\equiv c_1(x) $ (see also \cite{Smyshlyaev:2004} with $f(x,t)\not\equiv 0$), and for $  c(x,t)\equiv c_2(t)$,   respectively.

Let {$\lambda_0$} be an arbitrary  constant  satisfying
\begin{align}
\lambda_0>\sup_{(x,t)\in (0,1)\times  {\mathbb{R}_{>0}}} c(x,t) .\label{lambda0}
\end{align}
  Let $D:=\{(x,y)\in \mathbb{R}^2|~0\leq y \leq x\leq 1\}$ and
\begin{align}
 {\mu(x,y):=\lambda_0-c_1(x)+c_1(y)},\forall (x,y)\in D.\label{hatc}
 	\end{align}

In order to  design a boundary feedback law to stabilize  system \eqref{original system} in the sense of   different norms,  we
consider the kernel function $k: D\rightarrow \mathbb{R}$ that is independent of the time variable $t$ and satisfies
\begin{subequations}\label{kernel}
\begin{align}
 		  		  k_{xx}(x,y)-k_{yy}(x,y)
 	 = &  {\mu(x,y)}k(x,y)+f(x,y) 
  +\int_y^xk(x,z)f(z,y)\text{d}z,  \label{kb}\\
   2\frac{\text{d}}{\text{d}x}\left(k(x,x)\right)=&{\lambda_0},\label{ka} \\
 		  k_{y}(x,0)=&0,\label{kc}\\
   k(0,0)=&0,\label{kd}
 	\end{align}
\end{subequations}
 where   $\frac{\text{d}}{\text{d}x}\left(k(x,x)\right):=k_x(x,x)+k_y(x,x)$.  {The existence and regularity of  the kernel function $k$  that {satisfies} \eqref{kernel} will be {specified} in Section \ref{Sec. 4}.

Define the boundary feedback control law as
\begin{align}
 U(t):=-k(1,1)w(1,t)-\int_{0}^{1}k_x(1,y)w(y,t)\text{d}y.\label{control law}
\end{align}
Let $\beta\in(0,1)$ be a fixed constant. Define for the initial data the set
\begin{align*}
W_0:=\bigg\{&w|~w\in \mathcal{H}^{2+\beta}([0,1]),w_{x}(0)=0,  w_{x}(1)=-k(1,1)w(1)-\int_{0}^{1}k_x(1,y)w(y)\text{d}y \bigg\}.
\end{align*}
Concerning with the well-posedness and stability in different norms of system \eqref{original system}, we have the following theorem, which is
the main result obtained in this paper.
\begin{theorem}\label{theorem}
Given initial data in $  W_0$ and considering system \eqref{original system} under the feedback control law \eqref{control law},  the following statements hold true:
 \begin{enumerate}[(i)]
\item   system \eqref{original system}    admits a unique classical solution $w$, which belongs to $ \mathcal{H}^{2+\beta,1+\frac{\beta}{2}}([0,1]\times [0,T])$ and has the derivative $w_{xt}\in L^2((0,1)\times (0,T))$ for any $T\in \mathbb{R}_{>0}$;
 \item system \eqref{original system}   is exponentially stabilized in $L^p$-norm, and $W^{1,p}$-norm, having the {estimates} 
\begin{align}\label{(ii)-p1}
\|w[t]\|_{p}\leq C_1e^{-\underline{\lambda} t}\|w_0\|_{p},~\forall t\in {\mathbb{R}_{>0}},
\end{align}
{and}
\begin{align}\label{(ii)-p2}
\|w[t]\|_{{1,p}}\leq C_2e^{-\underline{\lambda} t}\|w_0\|_{{1,p}},~\forall t\in {\mathbb{R}_{>0}},
\end{align}
respectively,
where  $w_0\in W_0$ denotes the initial data,   $\underline{\lambda}$ is a positive constant given by
   \begin{align}\underline{\lambda}:=\inf_{(x,t)\in (0,1)\times {\mathbb{R}_{>0}}} (\lambda_0-c (x,t)),\label{underline-lambda}
\end{align} and $C_1,C_2$ are positive constants depending only on  $p,k$ and  the solution $l$ of \eqref{inversekernel}  for $p\in[1,+\infty)$, and on $k$ and $l$ for $p=+\infty$, respectively;

  \item the solution of system \eqref{original system}  is continuously dependent on    the $L^{p}$-norm,  and $W^{1,p}$-norm,  of the initial data,   having the  {estimates}
\begin{align*}
 &\|w_1-w_2\|_{C([0,T];L^p(0,1))} 
 \leq  C_1\| w_{01}-w_{02}\|_{p},\forall T\in \mathbb{R}_{>0},
 \end{align*}
and
 \begin{align*}
 &\|w_1-w_2\|_{C([0,T];W^{1,p}(0,1))} 
 \leq   C_2\| w_{01}-w_{02}\|_{1,p},\forall T\in \mathbb{R}_{>0},
 \end{align*}
 respectively,
    where $w_i  $ denotes the  solution of system \eqref{original system} with the  initial data $w_{0i}\in W_0  (i=1,2)$, and $C_1,C_2$  are the same as in (ii).
  \end{enumerate}

\end{theorem}


\begin{remark}\label{Remark 1}

 It is worth mentioning that under certain {Gevrey}-like conditions of the time-varying coefficient $c$, e.g., with the assumption that
$c[x]\in C^{\infty}( {\mathbb{R}_{>0}})$ and satisfies the following inequality for some  positive constant $R$ and  $i\in\mathbb{N}_0$:
 \begin{align*}
 \sup_{t\in  {\mathbb{R}_{>0}}}\left|\partial_t^i c(x,t)\right| \leq R^{i+1}(i!)^\alpha,\forall x\in [0,1],
 \end{align*}
the exponential  stabilization  in  $L^2$-norm was studied in, e.g., \cite{Kerschbaum:2020,meurer2009tracking,si2018boundary,smyshlyaev2005control}, for a class of time-varying systems with different boundary feedback controls. The kernel functions obtained in the existing literature depend on the time variable and   are differentiable (w.r.t. $t$) up to the order of $i$. Since $i$ can be given arbitrarily, the kernel functions  can be sufficiently smooth.

  In this paper, by using a combinatorial method, no Gevrey-like {condition} is  imposed on the time-varying coefficient $c$, and the stability is established in $L^p$-norm or $W^{1,p}$-norm for any $p\in[1,+\infty]$. In addition, the kernel functions considered in this paper are independent of the time variable (and therefore naturally infinitely differentiable  w.r.t. $t$), which  reduce   the associated computations  extensively.
\end{remark}

\subsection{Technical Line}\label{sec. 2.2}
 The proof of Theorem \ref{theorem} mainly consists of 3 steps that {are included} in Section \ref{Sec. 4}, Section \ref{Sec. 3} and Section \ref{Sec. 5}, respectively.

 {\textbf{Step 1:}  take into account an equivalent system  by using  transformations.}

 {Indeed, in order to avoid using a Gevrey-like condition, we transform  system \eqref{original system} into a target system that has also a space-time-varying coefficient. More precisely, {for   $\lambda_0$ satisfying \eqref{lambda0}, let first
 \begin{align}
 \lambda(x,t):=\lambda_0-c\left(x,t\right).\label{lambda}
 \end{align}
 Then using} the integral transformation
 \begin{align}\label{transformation}
 u(x,t):=w(x,t)+\int_{0}^{x}k(x,y)w(y,t)\text{d}y,
 \end{align}
 we  transform system \eqref{original system}   into the following target system:
\begin{subequations}\label{transformation system}
\begin{align}
&u_t(x,t)=u_{xx}(x,t)-\lambda (x,t)u(x,t),(x,t)\in (0,1)\times {\mathbb{R}_{>0}} ,\label{2a}\\
&u_x(0,t)=0,t\in {\mathbb{R}_{>0}}, \label{2b}\\
&u_x(1,t)=0, t\in {\mathbb{R}_{>0}},\label{2c}\\
&u(x,0)=u_0(x),x\in(0,1),
\end{align}
\end{subequations}
where $ u_0(x):=w_0(x)+\int_0^xk(x,y)w_0(y)\text{d}y$.

{Conversely, using the  integral transformation
\begin{align}\label{inverse}
w(x,t):=u(x,t)-\int_{0}^{x}l(x,y)u(y,t)\text{d}y,
\end{align}
we  transform system \eqref{transformation system} into the original system \eqref{original system}, where $l: D\rightarrow \mathbb{R}$  satisfies
\begin{subequations}\label{inversekernel}
\begin{align}
 		 		 l_{xx}(x,y)-l_{yy}(x,y)
 	 =&{\phi(x,y)}l(x,y)+f(x,y) -\int_y^xl(x,z)f(z,y)\text{d}z,\\
 2\frac{\text{d}}{\text{d}x}\left(l(x,x)\right)=&\lambda_0,\\
 		 l_{y}(x,0)=&0,\\
  l(0,0)=&0,
 	\end{align}
\end{subequations}}
  {with  $\phi(x,y):=-\lambda_0-c_1(x)+c_1(y)$.}

{The equivalence of {the original system} \eqref{original system} and system \eqref{transformation system} can be proved  in a direct way. {For instance, we} only show that if $w$ is the solution of   system \eqref{original system}, then $u$ is the solution of system \eqref{transformation system}. Indeed, the partial derivatives of (\ref{transformation}) w.r.t. $x$ and $t$ give
\begin{align*}
    u_x(0,t)=&w_x(0,t)+k(0,0)w(0,t),\\
     u_x(1,t)=&{w_x(1,t)+k(1,1)w(1,t)+\int_{0}^{1}k_x(1,y)w(y,t)\text{d}y}
  \end{align*}
  and}
  \begin{align*}
 &u_t(x,t)-u_{xx}(x,t)+\lambda(x,t) u(x,t)\\
 =&w(x,t)\bigg( -k_{y}(x,x)-\frac{\text{d}}{\text{d}x}\left(k(x,x)\right)-k_{x}(x,x) {+\lambda(x,t)+c(x,t)}\bigg)
 +\int_{0}^{x}w(y,t)\bigg(k_{yy}(x,y)-k_{xx}(x,y)\\
 & {+\lambda(x,t) k(x,y)+c(y,t)k(x,y)}
 +f(x,y) +\int_y^xk(x,z)f(z,y)\text{d}z\bigg)\text{d}y\\
 &
 +k_{y}(x,0)w(0,t) +\int_y^xk(x,z)f(z,y)\text{d}z\bigg)\text{d}y
 +k_{y}(x,0)w(0,t)\\
 =&w(x,t)\bigg(  -k_{y}(x,x)-\frac{\text{d}}{\text{d}x}\left(k(x,x)\right)-k_{x}(x,x)  {+\lambda_0}\bigg)  
 +\int_{0}^{x}w(y,t)\bigg(k_{yy}(x,y)-k_{xx}(x,y) \\
 & +\left(\lambda_0-c_1(x)-c_2(t)\right) k(x,y)+\left( c_1(y)+c_2(t)\right) k(x,y)
 +f(x,y)+\int_y^xk(x,z)f(z,y)\text{d}z\bigg)\text{d}y
 +k_{y}(x,0)w(0,t)\\
 =&w(x,t)\bigg(  -2\frac{\text{d}}{\text{d}x}\left(k(x,x)\right) {+\lambda_0}\bigg) 
 +\int_{0}^{x}w(y,t)\bigg(k_{yy}(x,y)-k_{xx}(x,y)  +\mu(x,y) k(x,y)
\\
 & +f(x,y)+\int_y^xk(x,z)f(z,y)\text{d}z\bigg)\text{d}y
 +k_{y}(x,0)w(0,t)\\
 =&0.
 \end{align*}
 Since $k$ satisfies \eqref{kernel}  and the boundary control law $U(t)$ is defined by \eqref{control law}, we deduce that $u$ satisfies \eqref{transformation system}.

{\textbf{Step 2:} establish   stability in $L^p$-norm, and  $W^{1,p}$-norm, for  the target system \eqref{transformation system} with initial data  belonging to $W_0$
  by using {the} {technique}  of {ALFs}.}

{As mentioned in Introduction,  singularities may appear  in the {Lyapunov} arguments when  $p\in[1,2)$. To overcome this obstacle, we  apply the {technique}  of {ALFs} that was   introduced in {\cite{zheng2021approximations,Zheng:202112}} for the analysis of integral input-to-state stability of nonlinear parabolic PDEs with external disturbances. Specifically,   for any $\tau\in  {\mathbb{R}_{>0}}$, let
\begin{align}\label{rho}
\rho_\tau(s):=\left\{\begin{aligned}
& |s|,\ |s|\geq \tau,  \\
& -\frac{s^4}{8\tau^3}+\frac{3s^2}{4\tau}+\frac{3\tau}{8},\ |s|< \tau.
\end{aligned}\right.
\end{align}
that is $C^2$-continuous in $s$ and satisfies for any $s\in \mathbb{R}$:
\begin{subequations}\label{property of rho}
	\begin{align}
	&\rho_\tau'(0)=0,~0\leq |s|\leq\rho_\tau(s),~|\rho_\tau'(s)|\leq 1,\label{property of rho a}\\
	&0\leq  \rho_\tau''(s)=\left\{\begin{aligned}
	& 0,\ |s|\geq \tau,\\
	&  \frac{3}{2\tau}\left(1-\frac{s^2}{\tau^2}\right),\ |s|< \tau,
	\end{aligned}\right. \label{property of rho b}\\
	&0\leq \rho_\tau(s)-\frac{3\tau}{8}\leq  \rho_\tau'(s)s\leq\rho_\tau(s) \leq |s|+\frac{3\tau}{8}. \label{property of rho c}
	\end{align}
\end{subequations}
The second inequality of (\ref{property of rho a}) and the fourth inequality of \eqref{property of rho c} guarantee that
\begin{align*}
&\int_{0}^{1}\rho_\tau^p(v)\text{d}x\rightarrow\int_{0}^{1}|v|^p\text{d}x  ~~\text{as}~~ \tau\rightarrow 0^+,\\
&\int_{0}^{1}\rho_\tau^p(v_x)\text{d}x\rightarrow\int_{0}^{1}|v_x|^p\text{d}x ~~\text{as}~~ \tau\rightarrow 0^+,
\end{align*}
for any $v\in W^{1,p}(0,1)$. Therefore, instead of  choosing  $\int_{0}^{1}|v|^p\text{d}x$ or $\int_{0}^{1}|v|^p\text{d}x+\int_{0}^{1}|v_x|^p\text{d}x$ as a Lyapunov functional, we take  into account the functional $\int_{0}^{1}\rho_\tau^p(v)\text{d}x$ or  $\int_{0}^{1}\rho_\tau^p(v)\text{d}x+\int_{0}^{1}\rho_\tau^p(v_x)\text{d}x$ to establish the stability in $L^p$-norm or  $W^{1,p}$-norm for the target system \eqref{transformation system}  in the Lyapunov arguments.}

{\textbf{Step 3:}  establish the exponentially stability  in $L^p$-norm, and  $W^{1,p}$-norm, for    the original system \eqref{original system} with initial data  in $W_0$ by using  the inverse transformation \eqref{inverse},  the  estimates of $u$, and the properties of the kernel functions {$k$ and} $l$.}

\section{Well-posedness and {Stability  of the Target System}}\label{Sec. 3}


{In accordance with $W_0$, we define
\begin{align*}
U_0:=\{u|~u\in \mathcal{H}^{2+\beta}([0,1]),u_{x}(0)=u_{x}(1)=0 \}
\end{align*}
for the initial data of the target system  \eqref{transformation system}.
\begin{proposition}\label{well-posedness}Consider system  \eqref{transformation system} with  initial data in $U_0$. The following statements hold true:
\begin{enumerate}[(i)]
 \item system  \eqref{transformation system}  admits a unique classical solution $u$, which belongs to
$ \mathcal{H}^{2+\beta,1+\frac{\beta}{2}}([0,1]\times [0,T])$ and  has the derivative $u_{xt}\in L^2((0,1)\times (0,T))$  for any $T>0$;
 \item  system \eqref{transformation system}   is exponentially stable in $L^p$-norm, and $W^{1,p}$-norm, {having} the   estimates
\begin{align}
\|u[t]\|_{p}\leq&e^{-\underline{\lambda}  t}\|u_0\|_{p},~\forall t\in {\mathbb{R}_{>0}},\label{stability1}
\end{align}
and
\begin{align}
\|u[t]\|_{{1,p}}\leq&e^{-\underline{\lambda}  t}\|u_0\|_{{1,p}},~\forall t\in {\mathbb{R}_{>0}},\label{stability2}
\end{align}
respectively,
 {where  where $u   $  denotes  the  solution of system \eqref{transformation system} corresponding to the initial data $u_{0}\in U_0$, $\underline{\lambda} $ is a positive constant defined by \eqref{underline-lambda};}
\item the solution of system  \eqref{transformation system}   depends continuously on  the $L^{p}$-norm,  and $W^{1,p}$-norm,  of the initial data, having the esimates 
    \begin{align}
  \|u_1-u_2\|_{C([0,T];L^p(0,1))} 
 \leq   \| u_{01}-u_{02}\|_{L^p(0,1)} ,\label{dependency}
 \end{align}
 and
 \begin{align}
 \|u_1-u_2\|_{C([0,T];W^{1,p}(0,1))} \leq   \| u_{01}-u_{02}\|_{W^{1,p}(0,1)} ,\label{dependency2}
 \end{align}
 respectively,
    where $u_i  $  denotes  the  solution of system \eqref{transformation system} corresponding to the initial data $u_{0i}\in U_0$, $i=1,2$.

\end{enumerate}
\end{proposition}}
\begin{pf*}{Proof.}
The assertion  (i)  follows from \cite[Chap. V, Theorem 7.4]{Ladyzhenskaya:1968} and the proof of  \cite[Chap. V, Lemma 7.2]{Ladyzhenskaya:1968}.


%

 We adopt the {technique of ALFs} to prove  the assertions (ii) and (iii).

  Indeed, let $u   $  be  the  solution of system \eqref{transformation system} with the initial data $u_{0}\in U_0$. For any $\tau\in  {\mathbb{R}_{>0}}$, let $\rho_{\tau}:\mathbb{R}\rightarrow \mathbb{R}_{>0}$ be defined by \eqref{rho}.

For $p\in [1,+\infty)$, by integrating by parts, we get
 \begin{align}
 \frac{1}{p}\frac{\text{d}}{\text{d}t}\int_0^1\rho^p_{\tau}(u)\text{d}x
 =&\int_0^1\rho^{p-1}_{\tau}(u)\rho_{\tau}'(u)u_t\text{d}x\notag\\
 =&\int_0^1\rho^{p-1}_{\tau}(u)\rho_{\tau}'(u)\left(u_{xx}-\lambda u\right)\text{d}x\notag\\
=&    (\rho^{p-1}_{\tau}(u)\rho_{\tau}'(u)u_x)|_{x=0}^{x=1}  - \int_0^1\left( (p-1) \rho_{\tau}^{p-2}(u)(\rho_{\tau}'(u))^2+\rho_{\tau}^{p-1}(u)\rho_{\tau}''(u)\right)u_x^2 \text{d}x      \notag\\&- \int_0^1\rho_{\tau}^{p-1}(u)\rho_{\tau}'(u)\lambda u\text{d}x.\label{total-int}
\end{align}

{Note} that
\begin{align}
\phi(s):=&  (p-1) \rho_{\tau}^{p-2}(s)(\rho_{\tau}'(s))^2+\rho^{p-1}_{\tau}(s)\rho_{\tau}''(s) 
 \geq   0,{\forall s\in \mathbb{R}.}\label{phi(s)}
\end{align}
{By \eqref{property of rho c}, \eqref{total-int} and \eqref{phi(s)}}, we deduce that
\begin{align*}
 \frac{1}{p}\frac{\text{d}}{\text{d}t}\int_0^1\rho_{\tau}^p(u)\text{d}x
\leq&- \int_0^1\rho_{\tau}^{p-1}(u)\rho_{\tau}'(u)\lambda u\text{d}x\notag\\
\leq&-\int_0^1{\lambda}\rho_{\tau}^{p-1}(u)\left(\rho_{\tau}(u)-\frac{3\tau}{8}\right)\text{d}x\notag\\
=&-\int_0^1{\lambda}\rho_{\tau}^{p}(u)\text{d}x+\frac{3\tau}{8}\int_0^1{\lambda}\rho_{\tau}^{p-1}(u)\text{d}x,
\end{align*}
{which {along with \eqref{lambda0}, \eqref{underline-lambda} and \eqref{lambda}} implies}
\begin{align}\label{pu}
 \frac{\text{d}}{\text{d}t}\int_0^1\rho_{\tau}^p(u)\text{d}x 
\leq&-\underline{\lambda} p\int_0^1\rho_{\tau}^{p}(u)\text{d}x+\frac{3}{8} \tau p\int_0^1{\lambda}\rho_{\tau}^{p-1}(u)\text{d}x.
\end{align}
By {Lemma \ref{lemma}} in Appendix, for   any $t\in [0,T]$, it holds that
 \begin{align} \label{u}
  \int_0^1 \rho_{\tau}^p (v(x,t))\text{d}x
 \leq&e^{-\underline{\lambda} pt}\int_0^1 \rho_{\tau}^{p} (u_0 (x))\text{d}x +\tau\int_0^t e^{-\underline{\lambda} p(t-s)}\psi_1(s)\text{d}s,
 \end{align}
where $\psi_1(s):=\frac{3}{8} p \int_0^1{\lambda}(x,s)\rho_{\tau}^{p-1}(u(x,s))\text{d}x$.

Letting $\tau \rightarrow 0^+$, we infer from  \eqref{u}   that
 \begin{align}\label{equ. 22}
\|u[t]\|_{p}\leq&e^{-\underline{\lambda}  t}\|u_0\|_{p},~\forall t \in [0,T],
\end{align}
which is \eqref{stability1} in the assertion (ii) for $p\in [1,+\infty)$.


  Analogously,  for $p\in [1,+\infty)$, we calculate
\begin{align}
 \frac{1}{p}\frac{\text{d}}{\text{d}t}\int_0^1\rho_{\tau}^{p}(u_x)\text{d}x 
=&\int_0^1\rho_{\tau}^{p-1}(u_x)\rho_{\tau}'(u_x){u_{xt}}\text{d}x\notag\\
=&\left(\rho_{\tau}^{p-1}(u_x)\rho_{\tau}'(u_x)u_{t}\right)|_{x=0}^{x=1}\ -\int_0^1\left(\rho_{\tau}^{p-1}(u_x)\rho_{\tau}'(u_x)\right)_{x}u_{t}\text{d}x
\notag\\
=&- \int_0^1 \left(\rho_{\tau}^{p-1}(u_x)\rho_{\tau}'(u_x)\right)_{x}  \left(u_{xx}-\lambda u\right)\text{d}x\notag\\
=&-\int_0^1\phi(u_x)u_{xx}^2 \text{d}x+
\int_0^1\left(\rho_{\tau}^{p-1}(u_x)\rho_{\tau}'(u_x)\right)_{x}\lambda u \text{d}x\notag\\
=&-\int_0^1\phi(u_x)u_{xx}^2 \text{d}x-\int_{0}^{1}\rho_{\tau}^{p-1}(u_x)\rho_{\tau}'(u_x)\lambda u_x\text{d}x.\label{total-int-2'}
\end{align}

{By \eqref{property of rho c}, \eqref{phi(s)} and \eqref{total-int-2'},} we get
\begin{align}\label{pux'}
&\frac{\text{d}}{\text{d}t}\int_0^1\rho_{\tau}^{p}(u_x)\text{d}x 
\leq -\underline{\lambda} p\int_0^1\rho_{\tau}^{p}(u_x)\text{d}x+\frac{3}{8}\tau p\int_0^1{\lambda}\rho_{\tau}^{p-1}(u_x)\text{d}x.
\end{align}
{By \eqref{pu} and \eqref{pux'}, we have}
 \begin{align}
  \frac{\text{d}}{\text{d}t}\left( \int_0^1\rho_{\tau}^{p}(u)\text{d}x+\int_0^1\rho_{\tau}^{p}(u_x)\text{d}x\right)  
 \leq&-\underline{\lambda} p\int_0^1\left(\rho_{\tau}^{p}(u)+\rho_{\tau}^{p}(u_x)\right)\text{d}x +\frac{3}{8}\tau p\int_0^1{\lambda}\left(\rho_{\tau}^{p-1}(u)+\rho_{\tau}^{p-1}(u_x)\right)\text{d}x.
 \end{align}
By {Lemma \ref{lemma}}  in Appendix, for  any $t\in[0,T]$, it holds that
 \begin{align}\label{ux'}
  \int_0^1 \left( \rho_{\tau}^{p} (u(x,t))+\rho_{\tau}^{p} (u_x(x,t))\right) \text{d}x 
 \leq e^{-\underline{\lambda} pt}\int_0^1\left( \rho_{\tau}^{p} (u_0 (x))+\rho_{\tau}^{p} (u_x(x,0))\right) \text{d}x +\tau \int_0^t e^{-\underline{\lambda} p(t-s)}\psi_2(s)\text{d}s,
 \end{align}
where $\psi_2(s):=\frac{3}{8} p \int_0^1{\lambda}\left( \rho_{\tau}^{p-1}(u(x,s))+\rho_{\tau}^{p-1}(u_x(x,s)\right) \text{d}x$.

 Letting $\tau \rightarrow 0^+$, we infer from   \eqref{ux'}   that
 \begin{align}\label{equ. 22'}
\|u[t]\|_{1,p}\leq&e^{-\underline{\lambda}  t}\|u_0\|_{1,p},~\forall t \in [0,T],
\end{align}
which is \eqref{stability2} in the assertion (ii) for $p\in [1,+\infty)$.

Now letting $p\rightarrow+\infty$, by the classical result  of functional analysis, e.g., \cite[pp.74]{F.Riesz:1990}, we deduce that \eqref{stability1} and  \eqref{stability2}  also hold true for    $p=+\infty$.

In order to prove the assertion (iii), letting $v:=u_1-u_2$ and $v_0:=u_{01}-u_{02}$, we infer from the linearity of system \eqref{transformation system} and the assertion (ii) that
\begin{align*}
\|v[t]\|_{p}\leq&e^{-\underline{\lambda}  t}\|v_0\|_{p},~\forall t\in {\mathbb{R}_{>0}},
\end{align*}
and
\begin{align*}
\|v[t]\|_{{1,p}}\leq&e^{-\underline{\lambda}  t}\|v_0\|_{{1,p}},~\forall t\in {\mathbb{R}_{>0}},
\end{align*}
hold  true for all $p\in [1,+\infty]$. It follows that
\begin{align*}
\max_{t\in[0,T]}\|v[t]\|_{p}\leq& \|v_0\|_{p},\forall T\in \mathbb{R}_{>0},
\end{align*}and
\begin{align*}
\max_{t\in[0,T]}\|v[t]\|_{1,p}\leq& \|v_0\|_{1,p},\forall T\in \mathbb{R}_{>0}.
\end{align*}
  Therefore, the assertion (iii) is true. The proof of Proposition \ref{well-posedness}  is complete.
$\hfill\blacksquare$

\end{pf*}

\section{{Well-posedness and Regularity of Kernel Functions}}\label{Sec. 4}
{In this {section}, we prove the existence, {uniqueness}  and regularity of kernel functions $k$ and $l$ by using the method of successive approximation {as in,  e.g., \cite{krstic2008boundary,smyshlyaev2005control,si2018boundary}.} Although the proof is standard, we would like to provide details for the  completeness. In particular, as indicated later in Remark \ref{Remark 2},  we fix a gap that exists in the proof of  uniqueness  of kernel functions in the existing literature.

\begin{proposition}\label{proposition kernel}{The problem  \eqref{kernel}, and \eqref{inversekernel}, has a unique solution $k$ and $l$, respectively, both of which belong to $\mathcal{H}^{2+\beta}\left(D\right)$.} 
\end{proposition}

\begin{pf*}{Proof.}
{Let}
$
\xi:=x+y$ and $\eta:=x-y.
$
{It follows that} $\eta\in[0,1]$ and $\xi\in[\eta,2-\eta]$. {Let}
\begin{align*}
G(\xi,\eta):=k(x,y)=k\left(\frac{\xi+\eta}{2},\frac{\xi-\eta}{2}\right).
\end{align*}
{Then} \eqref{kernel} is changed into
\begin{subequations}\label{Gset}
\begin{align}
4G_{\xi\eta}(\xi,\eta)
=& { \widetilde{\mu}(\xi,\eta)}G(\xi,\eta)+f \left(\frac{\xi+\eta}{2},\frac{\xi-\eta}{2}\right)
 +\int_{\frac{\xi-\eta}{2}}^{\frac{\xi+\eta}{2}}
 \widetilde{\mathcal{G}}(z,\xi,\eta)\text{d}z, \label{Gb}\\
 G(\xi,0)=&{\frac{\lambda_0}{4}}\xi, \label{Ga}\\
G_\xi(\xi,\xi)=&G_\eta(\xi,\xi), \label{Gc}\\
G(0,0)=&0, \label{Gd}
\end{align}
\end{subequations}
where
\begin{align*}
 &{\widetilde{\mu}(\xi,\eta):=\mu(x,y)=\lambda_0-c_1\left(\frac{\xi+\eta}{2}\right) +c_1\left(\frac{\xi-\eta}{2}\right)},\\ &\widetilde{\mathcal{G}}(z,\xi,\eta):=f\left(z,\frac{\xi-\eta}{2}\right)G\left(\frac{\xi+\eta}{2}+z,\frac{\xi+\eta}{2}-z\right).
\end{align*}
Integrating \eqref{Gb} w.r.t. $\eta$ from $0$ to $\eta$ and differentiating \eqref{Ga} w.r.t. $\xi$ give
\begin{align}\label{Ge}
G_\xi(\xi,\eta)=&{\frac{\lambda_0}{4}}+{{\frac{{ 1}}{4}}\int_{0}^{\eta}\widetilde{\mu}(\xi,s)G(\xi,s)\text{d}s} +\frac{1}{4}\int_{0}^{\eta}f\left(\frac{\xi+s}{2},\frac{\xi-s}{2}\right)\text{d}s +\frac{1}{4}\int_{0}^{\eta}\int_{\xi}^{\xi+\eta-s}\widehat{G}(\tau,s,\xi)\text{d}\tau\text{d}s,
\end{align}
where $\widehat{G}(\tau,s,\xi):=f\left(\frac{\tau-s}{2},\xi-\frac{\tau+s}{2}\right)G(\tau,s)$.

Integrating  {\eqref{Ge}} w.r.t. $\xi$ from $\eta$ to $\xi$ gives
\begin{align}\label{equ.411-3}
G(\xi,\eta)=&G(\eta,\eta)+{\frac{\lambda_0}{4}}(\xi-\eta) +{{\frac{{ 1}}{4}}\int_{\eta}^{\xi}\int_{0}^{\eta}\widetilde{\mu}(\tau,s)G(\tau,s)\text{d}s\text{d}\tau} +\frac{1}{4}\int_{\eta}^{\xi}\int_{0}^{\eta}f\left(\frac{\tau+s}{2},\frac{\tau-s}{2}\right)\text{d}s\text{d}\tau\notag\\
&+\frac{1}{4}\int_{\eta}^{\xi}\int_{0}^{\eta}\int_{z}^{z+\eta-s}\widehat{G}(\tau,s,z)\text{d}\tau\text{d}s\text{d}z.
\end{align}
By \eqref{Gc}, we have
\begin{align}
\frac{\text{d}}{\text{d}\xi}\left(G(\xi,\xi)\right)=2G_\xi(\xi,\xi).\label{equ.411-1}
\end{align}
  Substituting $\eta=\xi$ into \eqref{Ge}, {then integrating it over $[0,\xi]$ and applying \eqref{Gd} and \eqref{equ.411-1}, we obtain}
\begin{align*}
G(\xi,\xi)=&{\frac{\lambda_0}{2}}\xi+{{\frac{{ 1}}{2}}\int_{0}^{\xi}\int_{0}^{\tau}\widetilde{\mu}(\tau,s)G(\tau,s)\text{d}s\text{d}\tau} +\frac{1}{2}\int_{0}^{\xi}\int_{0}^{\tau}f\left(\frac{\tau+s}{2},\frac{\tau-s}{2}\right)\text{d}s\text{d}\tau\notag\\
&+\frac{1}{2}\int_{0}^{\xi}\int_{0}^{z}\int_{z}^{2z-s}\widehat{G}(\tau,s,z)\text{d}\tau\text{d}s\text{d}z,
\end{align*}
{which gives
\begin{align}\label{equ.411-2}
G(\eta,\eta)=&{\frac{\lambda_0}{2}}\eta+{{\frac{{ 1}}{2}}\int_{0}^{\eta}\int_{0}^{\tau}\widetilde{\mu}(\tau,s)G(\tau,s)\text{d}s\text{d}\tau} +\frac{1}{2}\int_{0}^{\eta}\int_{0}^{\tau}f\left(\frac{\tau+s}{2},\frac{\tau-s}{2}\right)\text{d}s\text{d}\tau\notag\\
&+\frac{1}{2}\int_{0}^{\eta}\int_{0}^{z}\int_{z}^{2z-s}\widehat{G}(\tau,s,z)\text{d}\tau\text{d}s\text{d}z.
\end{align}}
{We deduce by \eqref{equ.411-3} and \eqref{equ.411-2} that}
\begin{align}\label{otherkernel}
G(\xi,\eta)=&{\frac{\lambda_0}{4}}\left(\xi+\eta\right)
+\frac{1}{4}\int_{\eta}^{\xi}\int_{0}^{\eta}f\left(\frac{\tau+s}{2},\frac{\tau-s}{2}\right)\text{d}s\text{d}\tau +\frac{1}{2}\int_{0}^{\eta}\int_{0}^{\tau}f\left(\frac{\tau+s}{2},\frac{\tau-s}{2}\right)\text{d}s\text{d}\tau\notag\\
&+{{\frac{{ 1}}{4}}\int_{\eta}^{\xi}\int_{0}^{\eta}\widetilde{\mu}(\tau,s)G(\tau,s)\text{d}s\text{d}\tau} +{{\frac{{ 1}}{2}}\int_{0}^{\eta}\int_{0}^{\tau}\widetilde{\mu}(\tau,s)G(\tau,s)\text{d}s\text{d}\tau} +\frac{1}{4}\int_{\eta}^{\xi}\int_{0}^{\eta}\int_{z}^{z+\eta-s}\widehat{G}(\tau,s,z)\text{d}\tau\text{d}s\text{d}z\notag\\
&+\frac{1}{2}\int_{0}^{\eta}\int_{0}^{z}\int_{z}^{2z-s}\widehat{G}(\tau,s,z)\text{d}\tau\text{d}s\text{d}z.
\end{align}
We rewrite the integral equation \eqref{otherkernel}   as
\begin{align}\label{twoparts}
G(\xi,\eta)=G_0(\xi,\eta)+\Phi_G(\xi,\eta),
\end{align}
where  we defined
\begin{align*}
G_0(\xi,\eta):=&{\frac{\lambda_0}{4}}\left(\xi+\eta\right) +\frac{1}{4}\int_{\eta}^{\xi}\int_{0}^{\eta}f\left(\frac{\tau+s}{2},\frac{\tau-s}{2}\right)\text{d}s\text{d}\tau +\frac{1}{2}\int_{0}^{\eta}\int_{0}^{\tau}f\left(\frac{\tau+s}{2},\frac{\tau-s}{2}\right)\text{d}s\text{d}\tau,
\end{align*}
{and}  $\Phi_G(\xi,\eta)$   via
\begin{align}
{\Phi_H(\xi,\eta)} :=&{{\frac{{ 1}}{4}}\int_{\eta}^{\xi}\int_{0}^{\eta}\widetilde{\mu}(\tau,s)G(\tau,s)\text{d}s\text{d}\tau} +{{\frac{{ 1}}{2}}\int_{0}^{\eta}\int_{0}^{\tau}\widetilde{\mu}(\tau,s)G(\tau,s)\text{d}s\text{d}\tau} +\frac{1}{4}\int_{\eta}^{\xi}\int_{0}^{\eta}\int_{z}^{z+\eta-s}\widehat{H}(\tau,s,z)\text{d}\tau\text{d}s\text{d}z\notag\\
&+\frac{1}{2}\int_{0}^{\eta}\int_{0}^{z}\int_{z}^{2z-s}\widehat{H}(\tau,s,z)\text{d}\tau\text{d}s\text{d}z,\notag\\
\widehat{H}(\xi,\eta,z):=&f\left(\frac{\xi-\eta}{2},z-\frac{\xi+\eta}{2}\right)H(\xi,\eta),\label{hatH}
\end{align}
for  any  $\eta\in[0,1]$, $\xi\in[\eta,2-\eta]$, {$z\in[\eta,2-\eta]$}, and any function $H$.

Define the sequence  {$\{G_n(\xi,\eta)\}$} with $n\in \mathbb{N}_0$ via
\begin{align}\label{sequence}
G_{n+1}(\xi,\eta):=G_0(\xi,\eta)+{\Phi_{G_n}}(\xi,\eta),
\end{align}
and {let}
$\Delta G_{n}(\xi,\eta):=G_{n+1}(\xi,\eta)-G_{n}(\xi,\eta)$. It holds that
\begin{align*}
\Delta G_{n+1}(\xi,\eta)={\Phi_{\Delta G_n}}(\xi,\eta),
\end{align*}
and
\begin{align}\label{series sequence}
G_{n+1}{(\xi,\eta)}=G_0{(\xi,\eta)}+\sum_{j=0}^{n}\Delta G_{j}(\xi,\eta).
\end{align}
From \eqref{sequence}, we {see that if $\{G_n(\xi,\eta)\}$ converges uniformly w.r.t. $(\xi,\eta)$} when $n\rightarrow\infty$, then  $G(\xi,\eta):=\lim_{n\to\infty}G_n(\xi,\eta)$ is the solution of the integral equation {\eqref{twoparts}}, {and therefore the existence of a kernel function $k$ of \eqref{kernel} is guaranteed. In addition, in view of} \eqref{series sequence}, {the convergence of  $\{G_n(\xi,\eta)\}$ is equivalent to the convergence of the series $\sum_{n=0}^{\infty}\Delta G_n(\xi,\eta)$. So  it suffices} to show that the series $\sum_{n=0}^{\infty}\Delta G_n(\xi,\eta)$ is uniformly convergent {w.r.t. {$(\xi,\eta)$}. Furthermore, recalling  the Weierstrass M-test, it suffices to show that  the following inequality}
\begin{align}\label{induction}
\left|\Delta G_n(\xi,\eta)\right| \leq \frac{M^{n+2}}{(n+1)!}\left(\xi+\eta\right)^{n+1},
\end{align}
{holds true for all $n\in\mathbb{N}_0$,} where
\begin{align*}
 M:=&  {\frac{ \lambda_1+\overline{f}} {2}},  
\lambda_1:= 
 {\max\left\{|\lambda_0|,\max_{(x,y)\in[0,1]\times[0,1]}|\lambda_0-c_1(x)+c_1(y)| \right\}}, 
 {\overline{f}}:=  \max_{(x,y)\in[0,1]\times[0,1]}|f(x,y)|.
\end{align*}

 {Now we prove  \eqref{induction}  by induction. First of all, since}
\begin{align*}
\left| G_0(\xi,\eta)\right|&\leq{\frac{|\lambda_0|}{4}}\left(\xi+\eta\right)+\frac{{\overline{f}}}{4}(\xi\eta-\eta^2)+\frac{{\overline{f}}}{4}\eta^2\\
&\leq{\frac{|\lambda_0|}{4}}\left(\xi+\eta\right)+\frac{{\overline{f}}}{4}\xi\eta\\
&\leq{\frac{|\lambda_0|}{4}}\times2+\frac{{\overline{f}}}{4}\times2\times1\\
&\leq M,
\end{align*}
 {it follows that}
\begin{align}
\left|\Delta G_{0}(\xi,\eta)\right|=&\left|\Phi _{ G_{0}}(\xi,\eta) \right|\notag\\
\leq&\frac{{\lambda_1}}{4}\int_{\eta}^{\xi}\int_{0}^{\eta}M\text{d}s\text{d}\tau\notag +{\frac{{\lambda_1}}{2}}\int_{0}^{\eta}\int_{0}^{\tau}M\text{d}s\text{d}\tau +\frac{\overline{f}}{4}\int_{\eta}^{\xi}\int_{0}^{\eta}\int_{z}^{z+\eta-s}M\text{d}\tau\text{d}s\text{d}z +\frac{\overline{f}}{2}\int_{0}^{\eta}\int_{0}^{z}\int_{z}^{2z-s}M\text{d}\tau\text{d}s\text{d}z\notag\\
=&M\frac{{\lambda_1}}{4}(\xi\eta-\eta^2)+M\frac{{\lambda_1}}{4}\eta^2 +M\frac{\overline{f}}{8}(\xi\eta^2-\eta^3)+M\frac{\overline{f}}{12}\eta^3\notag\\
=&M\left(\frac{{\lambda_1}}{4}\xi\eta+\frac{{\overline{f}}}{8}\xi\eta^2-\frac{{\overline{f}}}{24}\eta^3\right)\notag\\
\leq&M\left(\frac{{\lambda_1}}{4}\xi\times1+\frac{{\overline{f}}}{8} {\times2 \times1 \times \eta -0} \right)\notag\\
\leq&{M\left(\frac{{\lambda_1}}{4} +\frac{\overline{f}}{4}\right)(\xi+\eta)}\notag\\
\leq&M^2\left(\xi+\eta\right),
\end{align}
which shows that \eqref{induction} holds for $n=0$.

{Supposing that the inequality \eqref{induction}  holds true for a general $n\in \mathbb{N}$, we need to show that \eqref{induction} also holds true for  $n+1$. Indeed,} we get
\begin{align} \label{en}
\left|\Delta G_{n+1}(\xi,\eta)\right|=&\left|{\Phi_{\Delta G_n}}(\xi,\eta)  \right| 
=  |I_1+I_2+I_3+I_4|,
\end{align}
where
\begin{align*}
I_1:=&{\frac{{ 1}}{4}}\int_{\eta}^{\xi}\int_{0}^{\eta}{\widetilde{\mu}(\tau,s)}\Delta G_{n}(\tau,s)\text{d}s\text{d}\tau,\notag\\
I_2:=&{\frac{{ 1}}{2}}\int_{0}^{\eta}\int_{0}^{\tau}{\widetilde{\mu}(\tau,s)}\Delta G_{n}(\tau,s)\text{d}s\text{d}\tau,\notag\\
I_3:=&\frac{1}{4}\int_{\eta}^{\xi}\int_{0}^{\eta}\int_{z}^{z+\eta-s}\widehat{\Delta G_{n}}(\tau,s,z)\text{d}\tau\text{d}s\text{d}z,\\
I_4:=&\frac{1}{2}\int_{0}^{\eta}\int_{0}^{z}\int_{z}^{2z-s}\widehat{\Delta G_{n}}(\tau,s,z)\text{d}\tau\text{d}s\text{d}z,
\end{align*}
 {with} $\widehat{\Delta G_{n}} $ defined via \eqref{hatH}.

{Since  \eqref{induction} holds true for $n$, it holds that}
\begin{align*}
|I_1|\leq&\frac{{\lambda_1}M^{n+2}}{4(n+1)!}\int_{\eta}^{\xi}\int_{0}^{\eta}(\tau+s)^{n+1}\text{d}s\text{d}\tau\notag\\
=&\frac{{\lambda_1}M^{n+2}}{4(n+1)!}\frac{1}{n+2}\int_{\eta}^{\xi}{\left((\tau+\eta)^{n+2}-\tau^{n+2}\right)}\text{d}\tau\notag\\
\leq&\frac{{\lambda_1}M^{n+2}}{4(n+1)!}\frac{1}{n+2}\int_{\eta}^{\xi}(\tau+\eta)^{n+2}\text{d}\tau\notag\\
=&\frac{{\lambda_1}M^{n+2}}{4(n+1)!}\frac{1}{(n+2)(n+3)}\left((\xi+\eta)^{n+3}-(2\eta)^{n+3}\right),\\
|I_2|\leq&\frac{{\lambda_1}M^{n+2}}{2(n+1)!}\int_{0}^{\eta}\int_{0}^{\tau}(\tau+s)^{n+1}\text{d}s\text{d}\tau\\
\leq&\frac{{\lambda_1}M^{n+2}}{2(n+1)!}\frac{1}{n+2}\int_{0}^{\eta}(2\tau)^{n+2}\text{d}\tau\\
=&\frac{{\lambda_1}M^{n+2}}{4(n+1)!}\frac{1}{(n+2)(n+3)}(2\eta)^{n+3},\\
|I_3|\leq&\frac{\overline{f} M^{n+2}}{4(n+1)!}\int_{\eta}^{\xi}\int_{0}^{\eta}\int_{z}^{z+\eta-s}(\tau+s)^{n+1}\text{d}\tau\text{d}s\text{d}z\\
\leq&\frac{\overline{f} M^{n+2}}{4(n+1)!}\frac{1}{n+2}\int_{\eta}^{\xi}\int_{0}^{\eta}(z+\eta)^{n+2}\text{d}s\text{d}z\\
=&\frac{\overline{f} M^{n+2}}{4(n+1)!}\frac{\eta}{n+2}\int_{\eta}^{\xi}(z+\eta)^{n+2}\text{d}s\text{d}z\\
\leq&\frac{\overline{f} M^{n+2}}{4(n+1)!}\frac{1}{n+2}\int_{\eta}^{\xi}(z+\eta)^{n+2}\text{d}z\\
=&\frac{\overline{f} M^{n+2}}{4(n+1)!}\frac{1}{(n+2)(n+3)}\left((\xi+\eta)^{n+3}-(2\eta)^{n+3}\right),\\
|I_4|\leq&\frac{\overline{f} M^{n+2}}{2(n+1)!}\int_{0}^{\eta}\int_{0}^{z}\int_{z}^{2z-s}(\tau+s)^{n+1}\text{d}\tau\text{d}s\text{d}z\\
\leq&\frac{\overline{f} M^{n+2}}{2(n+1)!}\frac{1}{n+2}\int_{0}^{\eta}\int_{0}^{z}(2z)^{n+2}\text{d}s\text{d}z\\
\leq&\frac{\overline{f} M^{n+2}}{2(n+1)!}\frac{2}{n+2}\int_{0}^{\eta}(2z)^{n+2}\text{d}z\\
=&\frac{\overline{f} M^{n+2}}{2(n+1)!}\frac{1}{(n+2)(n+3)}(2\eta)^{n+3}.
\end{align*}
Then we get
\begin{align*}
	 |\Delta G_{n+1}(\xi,\eta)| 
\leq &	\left(|I_1|+|I_2|\right)+\left(|I_3|+|I_4|\right)\\
\leq &\frac{{{\lambda_1}}}{4(n+3)}\frac{M^{n+2}}{(n+2)!}(\xi+\eta)^{n+3}  {+}\frac{{\overline{f}}}{2(n+3)}\frac{M^{n+2}}{(n+2)!}(\xi+\eta)^{n+3}\\
\leq & \left(\frac{{{\lambda_1}}}{2(n+3)}+\frac{{\overline{f}}}{n+3}\right)\frac{M^{n+2}}{(n+2)!}(\xi+\eta)^{n+2}\\
\leq &\frac{M^{n+3}}{(n+2)!}\left(\xi+\eta\right)^{n+2},
	\end{align*}
which implies that \eqref{induction} holds for $n+1$. Therefore, \eqref{induction} holds {true} for $n\in \mathbb{N}_0$.

 {We conclude} that the series
$
\sum_{n=0}^\infty \Delta G_n(\xi,\eta)
$
converges absolutely and uniformly in $\eta\in[0,1]$ and $\xi\in[\eta,2-\eta]$, and {$G(\xi,\eta):=\lim_{n\to\infty}G_n(\xi,\eta)$}  is a  solution of {\eqref{twoparts}. Moreover,  $G(\xi,\eta)$ is continuous in $\xi$ and $\eta$}. Furthermore, in view of \eqref{otherkernel} with $f\in C^1 \left([0,1]\times[0,1]\right)$ and $c_1\in C^1 \left([0,1] \right)$,   we deduce that $G_{\xi\xi},G_{\eta\eta},G_{\xi\eta}$ exist and are H\"{o}lder continuous with an arbitrary exponent $\beta'\in (0,1)$. In particular, $G_{\xi\xi},G_{\eta\eta},G_{\xi\eta}$ are H\"{o}lder continuous with the exponent $\beta$. Therefore, the problem  \eqref{kernel} admits a  solution   $k\in \mathcal{H}^{2+\beta}\left(D\right)$.


Regarding the  uniqueness of solution of \eqref{kernel}, which is equivalent to the uniqueness of  solution of the integral equation \eqref{otherkernel},  it suffices to show that any twice continuously differentiable solution  of  \eqref{otherkernel} can be approximated by the sequence $\{G_n(\xi,\eta)\}$ constructed via \eqref{sequence}. More precisely, for  any twice continuously differentiable solution   $\overline{G}(\xi,\eta)$ of \eqref{otherkernel} ,
    we intend   to show that the estimate
\begin{align} \label{unique}
|G_n(\xi,\eta)-\overline{G}(\xi,\eta)|\leq \frac{LM^{n}}{n!}\left(\xi+\eta\right)^{n},
\end{align}
holds true for all $n\in\mathbb{N}_0$, where  \begin{align*}
 L:= \max_{\eta\in[0,1],\xi\in[\eta,2-\eta]} |\Phi_{\overline{G}}(\xi,\eta) |.
 \end{align*}  This is because   \eqref{unique} implies the uniform  convergence of $\{G_n(\xi,\eta)\}$ to $ \overline{G}(\xi,\eta)$, while   the uniqueness of $ G(\xi,\eta):= \lim_{n\to\infty}G_n(\xi,\eta)$   leads to   $\overline{G}(\xi,\eta)\equiv G(\xi,\eta) $.    So the solution of \eqref{otherkernel} is
unique, and therefore the solution of \eqref{otherkernel} is also unique.

Now we prove \eqref{unique} by induction. First, since $\overline{G}(\xi,\eta)$ is a solution \eqref{otherkernel}, using the definition of $\Phi_{\overline{G}}$, we obtain
\begin{align*}
 |G_0(\xi,\eta)-\overline{G}(\xi,\eta)|=&|G_0(\xi,\eta)-\left(G_0(\xi,\eta)+ \Phi_{\overline{G}}(\xi,\eta)\right)| 
 = |\Phi_{\overline{G}}(\xi,\eta)| 
 \leq  L,
\end{align*}
which shows that \eqref{unique} holds true for $n=0$.

  Suppose that \eqref{unique} holds true for a general $n\in \mathbb{N}$. Then by  this assumption and the linearity of $\Phi_{H}(\xi,\eta)$  w.r.t. $H$, we have
\begin{align*}
|G_{n+1}(\xi,\eta)-\overline{G}(\xi,\eta)|= \left|     \Phi_{G_n}(\xi,\eta)  -\Phi_{\overline{G}}(\xi,\eta)\right | 
= \left|     \Phi_{G_n -\overline{G}}(\xi,\eta)\right | 
\leq I_5+I_6+I_7+I_8,
\end{align*}
where
\begin{align*}
I_5:=&L\frac{{\lambda_1}M^{n}}{4n!}\int_{\eta}^{\xi}\int_{0}^{\eta}(\tau+s)^{n}\text{d}s\text{d}\tau,\notag\\
I_6:=&L\frac{{\lambda_1}M^{n}}{2n!}\int_{0}^{\eta}\int_{0}^{\tau}(\tau+s)^{n}\text{d}s\text{d}\tau,\notag\\
I_7:=&L\frac{\overline{f} M^{n}}{4n!}\int_{\eta}^{\xi}\int_{0}^{\eta}\int_{z}^{z+\eta-s}(\tau+s)^{n}\text{d}\tau\text{d}s\text{d}z,\notag\\
I_8:=&L\frac{\overline{f} M^{n} }{2n!}\int_{0}^{\eta}\int_{0}^{z}\int_{z}^{2z-s}(\tau+s)^{n}\text{d}\tau\text{d}s\text{d}z.
\end{align*}
Analogous to the estimates of $I_1,I_2,I_3$ and $I_4$, it holds that
\begin{align*}
I_5 \leq&L\frac{{\lambda_1}M^{n}}{4n!}\frac{1}{(n+1)(n+2)}\left((\xi+\eta)^{n+2}-(2\eta)^{n+2}\right),\\
I_6 \leq&L\frac{{\lambda_1}M^{n}}{4n!}\frac{1}{(n+1)(n+2)}(2\eta)^{n+2},\\
I_7\leq&L\frac{\overline{f} M^{n}}{4n!}\frac{1}{(n+1)(n+2)}\left((\xi+\eta)^{n+2}-(2\eta)^{n+2}\right),\\
I_8 \leq&L\frac{\overline{f} M^{n} }{2n!}\frac{1}{(n+1)(n+2)}(2\eta)^{n+2}.
	\end{align*}
Then we obain
\begin{align*}
	 |G_{n+1}(\xi,\eta)-\overline{G}(\xi,\eta)| 
	\leq & L \left(\frac{{{\lambda_1}}}{2(n+2)}+\frac{{\overline{f}}}{n+2}\right)\frac{M^{n}}{(n+1)!}(\xi+\eta)^{n+1} 
\leq  {\frac{LM^{n+1}}{(n+1)!}\left(\xi+\eta\right)^{n+1}},
	\end{align*}
which implies that \eqref{unique} holds true for $n+1$. Thus,  \eqref{unique} holds true for all $n\in\mathbb{N}_0$. We conclude that the solution of \eqref{kernel} is unique.


Finally,  by  changing the variables and using the method of successive approximation as above, we may prove the well-posedness and regularity of the kernel function $l$ for \eqref{inversekernel}.$\hfill\blacksquare$
\end{pf*}
\begin{remark}\label{Remark 2}
For the stabilization problem of linear parabolic PDEs with time-independent coefficients, the ``uniqueness  of kernel functions'' was claimed in \cite{Smyshlyaev:2004} and proved by estimating
\begin{align}
 |\Delta G(\xi,\eta)|:= |G'(\xi,\eta)-G''(\xi,\eta)|\leq \frac{C}{n!}(\xi+\eta)^n\label{error}
  \end{align}
  with some constant $C>0$, where $G',G''$ are two kernel functions; see (39) and (40) of \cite{Smyshlyaev:2004}.

    It is worth noting that a gap exists in the proof of \cite{Smyshlyaev:2004}.
     On one hand,  if $G' $ and $G''$  are given by (see (39) of \cite{Smyshlyaev:2004})
\begin{align}\label{47}
 G'(\xi,\eta)=\sum_{n=0}^{\infty}G_n(\xi,\eta), ~
 G''(\xi,\eta)=\sum_{n=0}^{\infty}G_n(\xi,\eta),
  \end{align}
  as indicated in  \cite{Smyshlyaev:2004}, the estimate  \eqref{error} can be obtained in the same way as in (38) of \cite{Smyshlyaev:2004}.
   Therefore  the equality of  $G'=G''$ can be concluded by  letting $n\rightarrow \infty$. However, {since the limit of $\sum_{i=0}^{n}G_i(\xi,\eta)$ is unique, which is known by the limit theory, due to \eqref{47}}, the equality of $G'=G''$ holds true trivially. So the proof presented in \cite{Smyshlyaev:2004}
   is essentially for the uniqueness of the limit of $\{G_n\}$  rather than the uniqueness of solution.

        On the other hand, if $G'$ and $G''$ are not given by  \eqref{47}, {without imposing more conditions}, the estimate  \eqref{error} can not be obtained in the same way as in (38) of \cite{Smyshlyaev:2004}.

In order to  address the above-mentioned issue,     as presented in this paper, a feasible idea of proving the uniqueness of kernel functions  is to prove that all solutions of integral equations ({e.g., \eqref{otherkernel} or \eqref{twoparts})}   can be approximated by the sequence  $\{G_n\}$ constructed. {That is to say, $G'$ and $G''$ in \cite{Smyshlyaev:2004}    have to be with a form of \eqref{47}}.    Then  the uniqueness of solution {is determined} by the uniqueness of the limit of $\{G_n\}$. Note that,  by using this method, the estimate considered in this paper (see \eqref{unique}) is   different from  \eqref{error}. In particular, \eqref{unique} implies \eqref{error}.

\end{remark}

\begin{remark}\label{Remark 4.2}
 It is worth noting that for  general functions $f   $ and $c$ it is impossible to obtain  an analytic    solution $G$ of the integral equation \eqref{otherkernel} by using the iterating  formula\eqref{series sequence} (or \eqref{sequence}). In some special cases, the solution $G$ of \eqref{otherkernel} can be given by a refined iteration formula.  For example, consider $f\equiv0$ and $c_1(x)=rx^2$ with  a nonnegative constant $r$.  The integral equation \eqref{otherkernel} becomes
\begin{align}\label{new G}
G(\xi,\eta)= {\frac{\lambda_0}{4}}\left(\xi+\eta\right)
+{{\frac{{ 1}}{4}}\int_{\eta}^{\xi}\int_{0}^{\eta}(\lambda_0-r\tau s)G(\tau,s)\text{d}s\text{d}\tau} +{{\frac{{ 1}}{2}}\int_{0}^{\eta}\int_{0}^{\tau}(\lambda_0-r\tau s)G(\tau,s)\text{d}s\text{d}\tau},
\end{align}
whose solution is uniquely determined by
\begin{align}\label{48}
 G(\xi,\eta) 
=&\frac{\lambda_0}{4}(\xi+\eta)+\frac{\lambda_0}{4}\sum_{n=1}^{\infty}\sum_{i=0}^{n}\left(\frac{1}{4}\right)^{n}\lambda_0^iA_{n+1-i}^{n}T_{2n-i}(\xi,\eta),
\end{align}
where
\begin{align*}
{A_{n+1-i}^{n}}:=&\left\{
\begin{aligned}
&(-r)^n\Pi_{m=1}^{n}C_{2m}, && i=0,\\
&(A_{n-i}^{n-1}+A_{n+1-i}^{n-1})C_{2n-i}, && 0<i<n,\\
&\Pi_{m=1}^{n}C_{m}, && i=n,
\end{aligned}
\right.\\
C_{n}:=&\frac{1}{n(n+1)},\\
T_{2n-i}(\xi,\eta):=&\xi^{2n+1-i}\eta^{2n-i}+\xi^{2n-i}\eta^{2n+1-i}.
\end{align*}
The proof is provided  in Appendix.

\end{remark}

\section{Well-posedness and Stability   of the Original System}\label{Sec. 5}

%

In this section, we prove the well-posedness and exponential  stability   in different norms   for  the original system, i.e., system \eqref{original system}.

 {\paragraph*{Proof of Theorem \ref{theorem}.} The assertion of Theorem \ref{theorem} (i) follows immediately from the inverse transformation \eqref{inverse}, Proposition \ref{proposition kernel} and Proposition \ref{well-posedness}~(i).

For the proof of Theorem \ref{theorem} (ii), we define first the following qualities:}
\begin{align*}
  \alpha_1:=&\max\limits_{(x,y)\in D}|k(x,y)|,
 ~\beta_1 :=\max\limits_{(x,y)\in D}|l(x,y)|,\\
  \alpha_2:=&\max\limits_{x\in[0,1]}|k(x,x)|,
 ~\beta_2 :=\max\limits_{x\in[0,1]}|l(x,x)|,\\
  \alpha_3:=&\max\limits_{(x,y)\in D}|k_x(x,y)|,
 ~\beta_3:=\max\limits_{(x,y)\in D}|l_x(x,y)|.
\end{align*}
 {For $p\in[1,+\infty)$, we infer from \eqref{inverse} and the  Cauchy-Schwartz' inequality} that
 \begin{align}
\|w[t]\|_{p}^p\leq2^{p-1}(1+\beta_1^p) \|u[t]\|_{p}^p.\label{50}
\end{align}
Note that \begin{align*}
u_0(x) =w_0(x)+\int_0^xk(x,y)w_0(y)\text{d}y,
\end{align*}
 {which implies
\begin{align}
 \|u_0\|_{p}^p\leq&2^{p-1}(1+\alpha_1^p) \|w_0\|_{p}^p,\label{u0p}\\
  \|u_{0x}\|_{p}^p\leq&3^{p-1}(\alpha_2^p+\alpha_3^p)\|w_{0}\|_{p}^p+3^{p-1}\|w_{0x}\|_{p}^p.\label{u0xp}
 \end{align}
 By \eqref{stability1}, \eqref{50} and \eqref{u0p}, we obtain  the exponential stability   of   system \eqref{original system} in $L^p$-norm for $p\in[1,+\infty)$:
\begin{align}\label{C1}
\|w[t]\|_{p}\leq C_1e^{-\underline{\lambda} t}\|w_0\|_{p},~\forall t\in {\mathbb{R}_{>0}},
\end{align}
where $C_1:=\left(4^{p-1}(1+\alpha_1^p)(1+\beta_1^p)\right)^{\frac{1}{p}}$.}


{Now differentiating two-hand sides of \eqref{inverse} w.r.t. $x$, taking $L^p$-norm, and using the Cauchy-Schwartz' inequality, we have
\begin{align*}
\|w_x[t]\|_{p}^p\leq 3^{p-1}\|u_x[t]\|_{p}^p+3^{p-1}\left(\beta_2^p+\beta_3^p)\right\|u[t]\|_{p}^p,
\end{align*}
which along with \eqref{stability2}, \eqref{u0xp} and \eqref{C1}    implies the exponential stability   of   system \eqref{original system} in $W^{1,p}$-norm for $p\in[1,+\infty)$:
\begin{align}\label{C2}
	\|w[t]\|_{{1,p}}\leq C_2e^{-\underline{\lambda} t}\|w_0\|_{{1,p}},
\end{align}
where
$
 C_2:=\max\left\{9^{\frac{p-1}{p}}\gamma_1^{\frac{1}{p}},\gamma_2^{\frac{1}{p}}\right\}
$
} with
\begin{align*}
\gamma_1:=&\max\{1,\beta_2^p+\beta_3^p\},\\
 \gamma_2:=&C_1^p+ 9^{p-1}\gamma_1\left(1+\alpha_1^p+\alpha_2^p+\alpha_3^p\right).
\end{align*}
 Furthermore, letting $p\rightarrow+\infty$ in \eqref{C1} and \eqref{C2}, we obtain  the exponential stability  of   system \eqref{original system} in $L^\infty$-norm and $W^{1,\infty}$-norm, respectively:
 \begin{align*}
&\|w[t]\|_{\infty}\leq C_3e^{-\underline{\lambda} t}\|w_0\|_{\infty},~\forall t\in {\mathbb{R}_{>0}},\\
&\|w[t]\|_{{1,\infty}}\leq C_4e^{-\underline{\lambda} t}\|w_0\|_{{1,\infty}},~\forall t\in {\mathbb{R}_{>0}},
\end{align*}
where $C_3:=4{\gamma_3}$ and $C_4:=\max \{9{\gamma_4},C_1+9{\gamma_4}(1+\alpha_1+\alpha_2+\alpha_3) \}$
with
\begin{align*}
  {\gamma_3}:=&\left\{
\begin{aligned}
&1, && 0<\alpha_1\leq1,0<\beta_1\leq1,\\
&\beta_1, && 0<\alpha_1\leq1,\beta_1>1,\\
&\alpha_1, && \alpha_1>1,0<\beta_1\leq1,\\
&\alpha_1\beta_1, && \alpha_1>1,\beta_1>1,
\end{aligned}
\right.\\
{\gamma_4}:=&\max\{1,\beta_2+\beta_3\}.
\end{align*}
Therefore,  the assertion of Theorem \ref{theorem}~(ii) holds true.

Let $w_1  $ and $w_2$ be  the  solution of system \eqref{original system} with   initial data $w_{01}\in W_0$ and $w_{02}\in W_0$, respectively.  Using the linearity of system \eqref{original system} and Theorem \ref{theorem}~(ii),  the estimates \eqref{(ii)-p1} and \eqref{(ii)-p2} are satisfied with $w:=w_1-w_2$, which imply the  assertion of Theorem \ref{theorem}~(iii). The proof is {complete}.
$\hfill\blacksquare$

\section{Conclusion}\label{Sec. 6}

{In this paper, a combinatorial method, i.e., the method of backstepping and the technique of ALFs, was applied to address the exponential stabilization  and continuous dependence of solutions on initial
data for  a class of $1$-D space-time-varying linear parabolic PDEs  without imposing a   Gevrey-like condition on the time-varying coefficient. Compared to the results of the existing literature, the kernel functions obtained via the proposed method are independent of the time variable, therefore the associated calculations have been  extensively simplified. In addition, exponential stability and  continuous dependence of solutions on initial
data could be established
in different norms ($L^p$  and $W^{1,p}$ norms) by using the  {technique} of ALFs,  which was used to deal with singularities in the case of $p\in[1,2)$.}

{It is worth noting that, with a Gevrey-like condition, the proposed method is still suitable to  address the  exponential  stabilization and continuous dependence of solutions on initial
data in different norms for $1$-D linear parabolic PDEs  having a general form, e.g., $c(x,t)$ need not to be divided into two terms as     \eqref{c}. However,  without involving a   Gevrey-like condition, it is  still challenging to apply the method to obtain time-independent kernel functions and establish the stability estimates in different norms  for parabolic PDEs  having a general form, which will be studied in our future work.}

 \appendix
\section{Appendix}
\subsection*{A.1. The Bellman-Gronwall-Peano's  inequality}

{The following lemma used in this paper is concerned with    the Bellman-Gronwall-Peano's inequality, which can be found in, e.g., \cite[pp. 94]{Flett:1980}.}
\begin{lemma}\label{lemma}~~Let $T\in {\mathbb{R}_{>0}}$, if the function $z$ is a nonnegative, absolutely continuous function on $[0,T]$, and   satisfies the inequality
\begin{align*}
\frac{\text{d}z}{\text{d}t}\leq q(t)z(t)+h(t)  \ a.e.\ in \ [0,T],
\end{align*}
where $q$, $h\in L^1(0,T)$, then
\begin{align*}
z(t)\leq e^{\int_{0}^tq(s)\text{d}s}z(0)+\int_{0}^t  e^{ \int_{s}^tq(r)\text{d}r}h(s)\text{d}s,~\forall t\in[0,T].
\end{align*}
\end{lemma}
\subsection*{A.2. Proof of \eqref{50} in Remark \ref{Remark 4.2}}

  We  apply a new  approximation to obtain \eqref{50}, which is the solution of the integral equation \eqref{new G}.

Let $G_0(\xi,\eta):=\frac{\lambda_0}{4}\left(\xi+\eta\right)$, and
\begin{align*}
G_{n+1}(\xi,\eta):=&\frac{ 1}{4}\int_{\eta}^{\xi}\int_{0}^{\eta}(\lambda_0-r\tau s)G_{n}(\tau,s)\text{d}s\text{d}\tau +\frac{ 1}{2}\int_{0}^{\eta}\int_{0}^{\tau}(\lambda_0-r\tau s)G_{n}(\tau,s)\text{d}s\text{d}\tau.
\end{align*}
According to the  proof of Proposition \ref{proposition kernel}, the solution of \eqref{new G} exists and is unique. Therefore, as long as the series $\sum_{n=0}^{\infty}G_{n}(\xi,\eta)$ converges, $G(\xi,\eta):=\sum_{n=0}^{\infty}G_{n}(\xi,\eta)$ is the unique solution of \eqref{new G}.

%
%

 {Now we apply the induction to prove  the convergence of $\sum_{n=0}^{\infty}G_{n}(\xi,\eta)$ and the fact  that
\begin{align} \label{new Gn} G_{n}(\xi,\eta)=\frac{\lambda_0}{4}\left(\frac{1}{4}\right)^{n}\sum_{i=0}^{n}\lambda_0^iA_{n+1-i}^{n}T_{2n-i}(\xi,\eta).
\end{align}
Regarding the convergence of $\sum_{n=0}^{\infty}G_{n}(\xi,\eta)$, it suffices to show that
\begin{align}\label{new G induction}
|G_n(\xi,\eta)| \leq \frac{M_1^{n+1}}{(n+1)!}\left(\xi^{n+1}\eta^{n}+\xi^{n}\eta^{n+1}\right),\forall n\in\mathbb{N}_0,
\end{align}
  where $M_1:=r+|\lambda_0|$. Firstly, the following result can be established:
\begin{align*}
|G_0(\xi,\eta)|\leq&M_1\left(\xi+\eta\right),
\end{align*}
which shows that \eqref{new G induction} holds true for $n=0$.
The next step is to prove that if \eqref{new G induction} holds true for a general $n\in \mathbb{N}$, then \eqref{new G induction} also holds true for $n+1$. Indeed, by direct computations, we have
\begin{align*}
|G_{n+1}(\xi,\eta)|
\leq&\frac{ 1}{4}\int_{\eta}^{\xi}\int_{0}^{\eta}(|\lambda_0|+r\tau s)|G_{n}(\tau,s)|\text{d}s\text{d}\tau +\frac{ 1}{2}\int_{0}^{\eta}\int_{0}^{\tau}(|\lambda_0|+r\tau s)|G_{n}(\tau,s)|\text{d}s\text{d}\tau\\
\leq&\frac{1}{4}\frac{M_1^{n+1}}{(n+1)!}\frac{1}{n+2}
\left(\frac{|\lambda_0|}{n+1}+\frac{r\xi\eta}{n+3}\right) \left(\xi^{n+2}\eta^{n+1}+\xi^{n+1}\eta^{n+2}\right)\\
\leq&\frac{M_1^{n+2}}{(n+2)!}\left(\xi^{n+2}\eta^{n+1}+\xi^{n+1}\eta^{n+2}\right),
\end{align*}
which implies that \eqref{new G induction} holds for $n+1$. Thus, \eqref{new G induction} holds true for all $n\in \mathbb{N}_0$.  Therefore, $\sum_{n=0}^{\infty}G_{n}(\xi,\eta)$ converges, and $G:=\sum_{n=0}^{\infty}G_{n}(\xi,\eta)$ is the solution of \eqref{new G}.

We need to prove \eqref{new Gn}. Note that
\begin{align}\label{regularity}
 \int_{\eta}^{\xi}\int_{0}^{\eta}\left(\tau^{n+1}s^n+\tau^ns^{n+1}\right)\text{d}s\text{d}\tau+ 
 2\int_{0}^{\eta}\int_{0}^{\tau}\left(\tau^{n+1}s^n+\tau^ns^{n+1}\right)\text{d}s\text{d}\tau 
= \frac{1}{(n+1)(n+2)}\xi^{n+1}\eta^{n+1}(\xi+\eta).
\end{align}
For any function $\nu$,   define   $K$ via
\begin{align*}
K(\nu(\tau,s)):=&\int_{\eta}^{\xi}\int_{0}^{\eta}\nu(\tau,s)\text{d}s\text{d}\tau+
2\int_{0}^{\eta}\int_{0}^{\tau}\nu(\tau,s)\text{d}s\text{d}\tau.
\end{align*}
By \eqref{regularity}, we have
\begin{align}\label{K}
K(T_n(\tau,s))=C_{n+1}T_{n+1}(\xi,\eta),
\end{align}
where
\begin{align*}
&T_n(\tau,s):=\tau^{n+1}s^n+\tau^ns^{n+1},\\
&C_{n+1}:=\frac{1}{(n+1)(n+2)}.
\end{align*}
Then, it holds that
\begin{align}
&K(\lambda_0T_n(\tau,s))=\lambda_0C_{n+1}T_{n+1}(\xi,\eta), \label{linearK1}\\
&K(\tau sT_n(\tau,s))=C_{n+2}T_{n+2}(\xi,\eta).\label{linearK2}
\end{align}
By \eqref{K}, \eqref{linearK1} and \eqref{linearK2}, we obtain
\begin{align*}
G_{1}(\xi,\eta)=&\frac{ \lambda_0}{4}\frac{1}{4}K((\lambda_0-r\tau s)T_0(\tau,s))
=\frac{\lambda_0}{4}\frac{1}{4}\left(\lambda_0C_{1}T_{1}(\xi,\eta)
-rC_{2}T_{2}(\xi,\eta)\right),
\end{align*}
which implies that  \eqref{new Gn} holds true for $n=1$.

Supposing that  \eqref{new Gn} holds true for $1,2,\ldots,n$, it follows that
\begin{align*}
G_{n+1}(\xi,\eta)
=&\frac{1}{4}K\left((\lambda_0-r\tau s)\frac{\lambda_0}{4}\left(\frac{1}{4}\right)^{n}\sum_{i=0}^{n}\lambda_0^iA_{n+1-i}^{n}T_{2n-i}(\xi,\eta)\right)\\
=&\frac{\lambda_0}{4}\left(\frac{1}{4}\right)^{n+1}\left(\sum_{i=0}^{n}\lambda_0^{i+1}A_{n+1-i}^{n}C_{2n+1-i}T_{2n+1-i}(\xi,\eta)\right.\left.-\sum_{i=0}^{n}r\lambda_0^{i}A_{n+1-i}^{n}C_{2n+2-i}T_{2n+2-i}(\xi,\eta)\right)\\
=&\frac{\lambda_0}{4}\left(\frac{1}{4}\right)^{n+1}\left(\sum_{i=0}^{n-1}\lambda_0^{i+1}A_{n+1-i}^{n}C_{2n+1-i}T_{2n+1-i}(\xi,\eta)\right.\\
&+\lambda_0^{n+1}A_{1}^{n}C_{n+1}T_{n+1}(\xi,\eta)
-rA_{n+1}^{n}C_{2n+2}T_{2n+2}(\xi,\eta)\left.+\sum_{i=1}^{n}\lambda_0^{i}A_{n+1-i}^{n}C_{2n+2-i}T_{2n+2-i}(\xi,\eta)\right)\\
=&\frac{\lambda_0}{4}\left(\frac{1}{4}\right)^{n+1}\bigg(-rA_{n+1}^{n}C_{2n+2}T_{2n+2}(\xi,\eta)+\sum_{i=1}^{n}\lambda_0^{i}\left(A_{n+1-i}^{n}+A_{n+2-i}^{n}\right)C_{2n+2-i}T_{2n+2-i}(\xi,\eta)\\
&+\lambda_0^{n+1}A_{1}^{n}C_{n+1}T_{n+1}(\xi,\eta)\bigg)\\
=&\frac{\lambda_0}{4}\left(\frac{1}{4}\right)^{n+1}\sum_{i=0}^{n+1}\lambda_0^{i}A_{n+2-i}^{n+1}T_{2n+2-i}(\xi,\eta),
\end{align*}
which implies that \eqref{new Gn} holds true for $n+1$. Therefore, by   induction, \eqref{new Gn} holds true for all $n\geq1$.} The proof is complete.
}

\end{document}